\documentclass[12pt]{article}
\newtheorem{thm}{Theorem}[section]

\newtheorem{appl}[thm]{Proposition} 
\newtheorem{lemma}[thm]{Lemma}
\usepackage{amssymb, amsmath}

\begin{document}

\title{Magic Cayley-Sudoku Tables\thanks{Portions of this paper are from the first author's senior thesis at Western Oregon University.}}
\author{Rosanna Mersereau \\ Columbus, OH
\and
Michael B. Ward \\ Western Oregon University}
\date{}

\maketitle

\section{Introduction}

Inspired by the popularity of sudoku puzzles along with the well-known fact that the body of the Cayley table\footnote{That is, the operation table of the group without the borders.} of any finite group already has $2/3$ of the properties of a sudoku table in that each element appears exactly once in each row and in each column, Carmichael, Schloeman, and Ward [\ref{csw}] introduced \emph{Cayley-sudoku tables}.  A Cayley-sudoku table of a finite group $G$ is a Cayley table for $G$ the body of which is partitioned into uniformly sized rectangular blocks, in such a way that each group element appears exactly once in each block. For example, Table \ref{z9} is a Cayley-sudoku table for $\mathbb{Z}_9 := \{0,1,2,3,4,5,6,7,8\}$ under addition mod 9 and Table \ref{mcst} is a Cayley-sudoku table for $\mathbb{Z}_3 \times \mathbb{Z}_3$ where the operation is componentwise addition mod 3 (and ordered pairs $(a,b)$ are abbreviated $ab$).  In each case, we see that we have a Cayley table of the group partitioned into $3 \times 3$ blocks that contain each group element exactly once.

Lorch and Weld [\ref{LW}] defined a modular magic sudoku table as an ordinary sudoku table (with $0$ in place of the usual $9$) in which the row, column, diagonal, and antidiagonal sums in each $3 \times 3$ block in the table are zero mod 9. ``Magic'' refers, of course, to magic Latin squares which have a rich history dating to ancient times.  Table \ref{mmst} is an example.  Each $3 \times 3$ block is a Latin square by the sudoku conditions and one can check that each  block is a magic square when ordinary addition is replaced with addition mod 9.  For example, consider the second block in the top row of blocks.  The the first column sum is $7+8+3 = 0 \pmod 9$ and the antidiagonal sum is $6+0+3 = 0 \pmod 9$.  Similarly, all row, column, diagonal, and antidiagonal sums in each of the blocks are $0 \pmod 9$, making each of them magic.

This paper aims to wed these two notions in a generalized situation.  That is, we study Cayley-sudoku tables of groups (not just groups of order 9) in which the blocks are magic squares with respect to {``adding'' the rows, columns, diagonals, and antidiagonals using the group operation}.  Table \ref{mcst} illustrates the idea for the group $\mathbb{Z}_3 \times \mathbb{Z}_3$.  All the relevant sums in each $3 \times 3$ block are $00$.  In fact, there is more.  Each block is a pandiagonal magic square, meaning the ``broken" diagonals and antidiagonals also sum to $00$.  Consider, for instance, the upper leftmost block.  Beginning in the second row, first column, the broken diagonal sum is $01 + 12 + 20 = 00$.

Although every finite group has at least one Cayley-sudoku table (see Construction 1 of [\ref{csw}]), not every finite group admits a magic Cayley-sudoku table.  First of all, in order to partition the Cayley table into square blocks, the order of the group must be a perfect square.  However, that alone is not sufficient.  We will show that there is no magic Cayley-sudoku table for $\mathbb{Z}_9$, which explains why we switched to $\mathbb{Z}_3 \times \mathbb{Z}_3$ in order to give an example.  We do not know a characterization of groups having magic Cayley-sudoku tables.  We are, however, able to give sufficient conditions which we apply to certain abelian groups, certain groups of prime power order, and which also allow us to show every finite group is isomorphic to a subgroup of a finite group admitting a magic Cayley-sudoku table.

\begin{table}[h!]
\begin{center}
\begin{tabular}{c||ccc|ccc|ccc|}
\null & 0 & 3 & 6 & 1 & 4 & 7 & 2 & 5 & 8 \\
\hline\hline 0 & 0 & 3 & 6 & 1 & 4 & 7 & 2 & 5 & 8 \\
1 & 1 & 4 & 7 & 2 & 5 & 8 & 3 & 6 & 0 \\
2 & 2 & 5 & 8 & 3 & 6 & 0 & 4 & 7 & 1 \\
\hline 3 & 3 & 6 & 0 & 4 & 7 & 1 & 5 & 8 & 2 \\
4 & 4 & 7 & 1 & 5 & 8 & 2 & 6 & 0 & 3 \\
5 & 5 & 8 & 2 & 6 & 0 & 3 & 7 & 1 & 4 \\
\hline 6 & 6 & 0 & 3 & 7 & 1 & 4 & 8 & 2 & 5 \\
7 & 7 & 1 & 4 & 8 & 2 & 5 & 0 & 3 & 6 \\
8 & 8 & 2 & 5 & 0 & 3 & 6 & 1 & 4 & 7 \\ \hline
\end{tabular}
\end{center}
\caption{$\mathbb{Z}_9$ Cayley-Sudoku Table} \label{z9}
\end{table}

\begin{table}[h!]
\begin{center}
\begin{tabular}{|ccc|ccc|ccc|} \hline
  1 & 8 & 0 & 7 & 5 & 6 & 4 & 2 & 3 \\
  2 & 3& 4&8&0&1&5&6&7 \\
  6&7&5&3&4&2&0&1&8 \\ \hline
  8&4&6&5&1&3&2&7&0 \\
  7&0&2&4&6&8&1&3&5 \\
  3&5&1&0&2&7&6&8&4 \\ \hline
  5&1&3&2&7&0&8&4&6 \\
  4&6&8&1&3&5&7&0&2 \\
  0&2&7&6&8&4&3&5&1 \\ \hline
\end{tabular}
\end{center}
\caption{Modular Magic Sudoku Table} \label{mmst}
\end{table}

\begin{table}
\begin{small} \begin{center}
\begin{tabular}{c||ccc|ccc|ccc|}
\null & 00 & 10 & 20 & 01 & 11 & 21 & 02 & 12 & 22 \\ \hline \hline
  00 & 00 & 10 & 20 & 01 & 11 & 21 & 02 & 12 & 22 \\
  01 & 01 & 11 & 21 & 02 & 12 & 22 & 00 & 10 & 20 \\
  02 & 02 & 12 & 22 & 00 & 10 & 20 & 01 & 11 & 21 \\ \hline
  10 & 10 & 20 & 00 & 11 & 21 & 01 & 12 & 22 & 02 \\
  11 & 11 & 21 & 01 & 12 & 22 & 02 & 10 & 20 & 00 \\
  12 & 12 & 22 & 02 & 10 & 20 & 00 & 11 & 21 & 01 \\ \hline
  20 & 20 & 00 & 10 & 21 & 01 & 11 & 22 & 02 & 12 \\
  21 & 21 & 01 & 11 & 22 & 02 & 12 & 20 & 00 & 10 \\
  22 & 22 & 02 & 12 & 20 & 00 & 10 & 21 & 01 & 11 \\ \hline
  \hline
\end{tabular}
\end{center} \end{small}
\caption{Magic Cayley-Sudoku Table for $\mathbb{Z}_3 \times \mathbb{Z}_3$} \label{mcst}
\end{table}

\section{Notation and Terminology}

 All groups considered herein are finite.  Multiplicative notation is used for abstract groups, with $1$ denoting the both the identity and the trivial subgroup $\{1\}$.  For abstract group elements $a$ and $b$, $a^b$ denotes $b^{-1}ab$.  When $H$ is a subgroup of $G$, we write $H \leq G$.  Ordered lists of elements are indicated with square brackets, $[a,b,c, \ldots]$. Furthermore, if $L = [a_1, a_2, \ldots, a_k]$ is a list of elements in an abstract group and $b$ is a group element, then $[Lb]$ denotes the list $[a_1 b, a_2 b, \ldots, a_k b]$.  On the other hand, when $H \leq G$, $[G:H]$ equals the index of $H$ of $G$. The \emph{exponent} of $G$, $\textrm{exp}(G)$, is the smallest positive integer such that $g^n = 1$ for all $g \in G$.  Finally, $Z(G)$ denotes the center of $G$.
 
\subsection*{Definitions and Background}

\begin{itemize}

\item For $H \leq G$ with $[G:H] = n$, if $g_1 H, g_2 H, \ldots g_n H$ are all the distinct left cosets of $H$ in $G$, we say $g_1, g_2, \ldots g_n$ is a \emph{complete set of left coset representatives}\footnote{The term``left transversal'' is also used.  We avoid it because it has a different meaning in the theory of Latin squares.} of $H$ in $G$.  In other words, $g_1, g_2, \ldots g_n$ is a collection of elements, exactly one from each left coset of $H$ in $G$.   A  \emph{complete set of right coset representatives} is defined analogously.

\item In any $k \times k$ array, a \emph{broken diagonal} is a list of entries in positions $(\ell,\ell+j)$ where $j$ is a fixed integer between $1$ and $k$; $1 \leq \ell \leq k$; and addition is mod $k$, writing $k$ in case of congruence to $0$ mod $k$.  Similarly, a \emph{broken antidiagonal} is a list of entries in positions $(\ell,j - \ell)$, with the same conditions on $j$ and $\ell$ and arithmetic mod $k$.
\bigskip

Informally, broken diagonals are parallel to the main diagonal of the square, but they wrap around when reaching a border of the square.  Broken antidiagonals behave the same way with respect to the main antidiagonal.  Note that $j=k$ gives the main diagonal, while $j=1$ gives the main antidiagonal.  Therefore, references to broken (anti)diagonals are understood to include main (anti)diagonals as well.

\item Operating elements from left to right in a row of an array of group elements is called a \emph{row product} of the array (or \emph{row sum} if the group is written additively).

\item Operating elements from top to bottom in a column, broken diagonal, or broken antidiagonal of an array of group elements is called a \emph{column, broken diagonal, or broken antidiagonal product}, respectively, of the array (\emph{sum} in an additive group).

\item  A \emph{Cayley-sudoku table} of a finite group $G$ is a Cayley table for $G$, the body of which is partitioned into uniformly sized rectangular blocks in such a way that each group element appears exactly once in each block.  

\item A Cayley-sudoku table is \emph{pandiagonal magic} provided the blocks are square and every row, column, broken diagonal, and broken antidiagonal product of each of the blocks gives the group identity.  

\end{itemize}

Before proceeding,we adopt a convenient convention for labeling a Cayley table.  When a list appears in a row or column of the table, it is to be interpreted as the individual elements of that list appearing (as ordered in the list)  in separate rows or columns, respectively.  For example, under that convention, the rows and columns of Table \ref{z9} could be labeled
\begin{center}\begin{tabular}{c||c|c|c|}
  & $[0,3,6]$ & $[1,4,7]$ & $\ldots$  \\
\hline \hline $[0,1,2]$ &  &  &  \\
\hline $[3,4,5]$ &  &  &  \\
\hline $\vdots$ &  &  &  \\
\hline
\end{tabular}. \end{center}
Thus, the label $[0,1,2]$ is interpreted as the elements 0, 1, and 2 listed vertically, one per row in the indicated order, and $[0,3,6]$ is interpreted as the elements 0, 3, and 6 listed horizontally, one per column in the indicated order, and so forth. 

For an introduction to Cayley-sudoku tables, see [\ref{csw}], from which we quote one construction to be used herein.

\begin{thm}[{Construction 1 of [\ref{csw}]} ] \label{const1} Let $G$ be a finite group. Assume $H$ is a subgroup of $G$ having order $k$ and index $n$ (so that $|G| = nk$).  Further suppose $Hg_{1}, Hg_{2},\ldots,Hg_{n}$ are the $n$ distinct right cosets of $H$ in $G$ and  $T_{1}, T_{2}, \ldots, T_{k}$ partition G into complete sets of left coset representatives of $H$ in $G$.  Then, for arbitrary orderings of the elements in each right coset and and in each complete set of left coset representatives, arranging the Cayley table of $G$ with columns labeled by the cosets $[Hg_{1}], [Hg_{2}],\ldots,[Hg_{n}]$ and the rows labeled by sets $[T_{1}], [T_{2}], \ldots, [T_{k}]$ (as in the table below) yields a Cayley-sudoku table of $G$ with blocks of dimension $n\times k$.

\begin{table}[!ht]
\begin{center}\begin{tabular}{c||c|c|c|c|}
  & $[Hg_{1}]$ & $[Hg_{2}]$ & $\ldots$ & $[Hg_{n}]$ \\
\hline \hline $[T_{1}]$ &  &  &  &  \\
\hline $[T_{2}]$ &  &  &  &  \\
\hline $\vdots$ &  &  &  &  \\
\hline $[T_{k}]$ &  &  &  &  \\
\hline
\end{tabular} \end{center}
\end{table}

\end{thm}

\section{The Principal Construction} \label{pconst}

As we shall see, the hypotheses in the following theorem are tailor made to produce the desired magic.

\begin{thm}[Magic Cayley-Sudoku Construction] \label{mcs}
Suppose

\begin{enumerate}
\item \label{exp} $G$ is a group, $|G| = k^2$, $\mathrm{exp}(G)$ divides $k$;
\item $N \leq Z(G)$, $|N| = k$;
\item \label{tpp} $\prod_{\ell =1}^k n_\ell = 1$ where $N = [n_1, n_2, \ldots n_k]$ (the order is unimportant since $N$ is abelian); and  
\item \label{msp} there is a complete set of left coset representatives $T := [t_1, t_2, \ldots t_k]$ of $N$ in $G$ such that for every $i$, $1 \leq i \leq k$, $\displaystyle{\prod_{\ell =1}^k (t_\ell t_i) = t_1 t_i t_2 t_i \cdots t_k t_i = 1}$.
\end{enumerate}

Then the following is a pandiagonal magic Cayley-sudoku table.

\begin{center}\begin{tabular}{c||c|c|c|c|}
  & $[N t_{1}]$ & $[N t_{2}]$ & $\ldots$ & $[N t_{k}]$ \\
\hline \hline $[Tn_{1}]$ &  &  &  &  \\
\hline $[Tn_{2}]$ &  &  &  &  \\
\hline $\vdots$ &  &  &  &  \\
\hline $[Tn_{k}]$ &  &  &  &  \\
\hline
\end{tabular}

\bigskip
Table M
\end{center}

\end{thm}

Henceforth, we refer to Table M as a \emph{table of $G$ based on $N$ and $T$}.
\bigskip

\emph{Proof:}  It is not hard to show $Tn_1, Tn_2, \ldots, Tn_k$ partition $G$ into complete sets of left coset representatives of $N$ in $G$.  Moreover, since $N$ normal subgroup of $G$, $T$ is also a complete set of right coset representatives of $N$ in $G$.  By Theorem \ref{const1}, Table M is a Cayley-sudoku table.  It remains to show Table M is pandiagonal magic.

The expansion of a typical block in Table M is

\begin{tabular}{c||c|}
  & $[N t_{i}]$  \\
\hline \hline $[Tn_m]$ & \null \\
\hline
\end{tabular}  =
\begin{tabular}{c||c|c|c|c|}
  & $n_1 t_{i}$ & $n_2 t_{i}$ & $\ldots$ & $n_k t_{i}$ \\
\hline \hline $t_{1} n_m$ &  &  &  &  \\
\hline $t_{2}n_m$ &  &  &  &  \\
\hline $\vdots$ &  &  &  &  \\
\hline $t_{k}n_m$ &  &  &  &  \\
\hline
\end{tabular}. 

Fix $j$ with $1 \leq j \leq n$.  

Since $N \leq Z(G)$, the $j^{th}$ row product is $\displaystyle{\prod_{\ell =1}^k (t_j n_m)(n_\ell t_i) = n_m^k (\prod_{\ell =1}^k n_\ell) (t_j t_i)^k }$. Furthermore, $n_m^k = 1 = (t_j t_i)^k$ and $\prod_{\ell =1}^k n_\ell = 1$ by hypotheses \ref{exp} and \ref{tpp}, respectively.  Therefore, the $j^{th}$ row product is $1$ as was to be shown.

Similarly, the $j^{th}$ column product is $$\displaystyle{\prod_{\ell =1}^k (t_\ell n_m )(n_j t_i) = n_m^k n_j^k \prod_{\ell =1}^k (t_\ell t_i)  = \prod_{\ell =1}^k (t_\ell t_i)  = 1 }$$ with the last equality coming from hypothesis \ref{msp}.

The broken diagonal product corresponding to $j$ is $$\displaystyle{\prod_{\ell =1}^k (t_\ell n_m )(n_{j+\ell} t_i) = n_m^k (\prod_{\ell =1}^k n_{j+\ell})(\prod_{\ell =1}^k (t_\ell t_i) )} = 1$$ where the middle product is $1$ by hypothesis \ref{tpp} with the factors cyclically permuted (or because $N$ is abelian) and, as before, the other factors are $1$ by hypotheses \ref{exp} and \ref{msp}. Likewise, the corresponding broken antidiagonal product is $$\displaystyle{\prod_{\ell =1}^k (t_\ell n_m )(n_{j-\ell} t_i) = n_m^k (\prod_{\ell =1}^k n_{j-\ell})(\prod_{\ell =1}^k (t_\ell t_i) ) = 1}$$.

(Broken diagonal and antidiagonal products were defined to be top to bottom.  However, in this particular situation, calculating left to right gives the same result.  For the left-to-right broken diagonal products we have $\displaystyle{\prod_{\ell =1}^k (t_{\ell-j} n_m)(n_\ell t_i) = n_m^k (\prod_{\ell =1}^k n_\ell)(\prod_{\ell =1}^k (t_{\ell-j} t_i) ) = 1  }$ with the last equality following from hypothesis \ref{msp} with factors cyclically permuted.  Left-to-right broken antidiagonal products are trivial in the same way.) 
\hspace{\fill} $\Box$
\bigskip

\section{Abelian Groups}

Here we investigate the existence and nonexistence of magic Cayley-sudoku tables for certain abelian groups.

\begin{thm}[Abelian Groups] \label{ag}  Suppose

\begin{enumerate}
\item \label{expa} $G$ is an abelian group, $|G| = k^2$, $\mathrm{exp}(G)$ divides $k$;
\item $N \leq G$, $|N| = k$;
\item \label{tppa} the product of the elements of $N$ is $1$; and
\item \label{mspa} the product of the elements of $G/N$ is $N$.
\end{enumerate}

Then there is a complete set of left coset representatives, $T$, of $N$ in $G$ such that the table of $G$ based on $N$ and $T$ is a pandiagonal magic Cayley-sudoku table.
\end{thm}

\emph{Proof.}  Consider any complete set of left coset representatives of $N$ in $G$, say $\{t_1', t_2', \ldots , t_k'\}$.  By hypothesis \ref{mspa}, $\displaystyle{\prod_{\ell =1}^k  t_\ell' N = N}$, which implies there exists $n \in N$ such that $\displaystyle{(\prod_{\ell =1}^k  t_\ell') n = 1}$.  Let $t_\ell = t_\ell'$ for $1 \leq \ell \leq k-1$ and $t_k = t_k'n$.  Clearly, $T := [t_1, t_2, \ldots , t_k]$ is a complete set of left coset representatives of $N$ in $G$ and $\displaystyle{\prod_{\ell =1}^k  t_\ell = 1}$.  We claim $T$ satisfies \ref{msp} of Theorem \ref{mcs}.  For any $i$, $1 \leq i \leq k$, $\displaystyle{\prod_{\ell =1}^k  (t_\ell t_i) = (\prod_{\ell =1}^k  t_\ell) t_i^k = 1 }$ because $G$ is abelian of exponent dividing $k$.  The rest now follows from Theorem \ref{mcs}.
\hspace{\fill} $\Box$

\bigskip
The following lemma shows when the hypotheses on the product of elements of $N$ and of $G/N$ in Theorem \ref{ag} are satisfied.  We will also use it in the next section.

\begin{lemma} \label{s2} Suppose $A:= [g_1, g_2, \ldots , g_n]$ is an abelian finite group.  Then $\displaystyle{\prod_{\ell =1}^n g_\ell = 1 }$ if and only if $A$ either has no elements of order $2$ or more than one element of order $2$ (equivalently, the Sylow 2-subgroups of $A$ are trivial or non-cyclic).
\end{lemma}

\emph{Proof.}  The proof is elementary; see Lemma 1 of [\ref{P}], for instance.
\hspace{\fill} $\Box$
\bigskip

In view of Lemma \ref{s2}, Theorem \ref{ag} may be applied, for example, to $G = \mathbb{Z}_9 \times \mathbb{Z}_3 \times \mathbb{Z}_3 \times \mathbb{Z}_4 \times \mathbb{Z}_4$ with $N = <3> \times \mathbb{Z}_3 \times \{0\} \times <2> \times <2>$.
\bigskip

We conclude this section by showing neither Theorem \ref{ag} nor any other method can produce a magic Cayley-sudoku table for some groups.   We use those groups to illustrate the necessity of some hypotheses in Theorem \ref{ag}.  


\begin{lemma} \label{nomcs}
There is no magic Cayley-sudoku table for $\mathbb{Z}_9$ nor for  $\mathbb{Z}_2 \times  \mathbb{Z}_2$.
\end{lemma}

\emph{Proof:}  Suppose to the contrary that there is a magic Cayley-sudoku table for $\mathbb{Z}_9$.  The magic blocks must be $3 \times 3$.  Consider the magic block intersecting the row and the column headed by $0$ in the table. Let $a$ and $b$ head the remaining two columns for that block and let $x$ and $y$ label the remaining two rows. (For example, the magic block might be as below.) By the magic property, the column sums are equal.  Hence, regardless of the order of the row and column labels, $x+y = x+y+3a = x+y+3b$, which implies $3a = 3b = 0$.  Thus, $a, b \in <3>$.  Similarly, equating the row sums yields $x, y \in <3>$.  However, this implies all the elements of the block are in $<3>$ contradicting the sudoku condition that each block must contain all the elements of $\mathbb{Z}_9$.

\begin{table}[h!]
\begin{center}
$\begin{array}{c||c|ccc|c|}
\null & \cdots & a & 0 & b & \cdots \\ \hline \hline
\vdots & \null & \null & \vdots & \null & \null \\ \hline
x & \null & x+a & x & x+b &  \null \\
y & \cdots & y+a & y & y+b & \cdots \\
0 & \null & a & 0 & b & \null \\ \hline
\vdots & \null & \null & \vdots & \null & \null \\
\end{array}$
\end{center}
\end{table}

Similarly, assume there is a magic Cayley-sudoku table for $\mathbb{Z}_2 \times \mathbb{Z}_2$.  Here, the magic blocks must be $2 \times 2$.  Consider the block intersecting the row and the column headed by $00$ and let $a$ and $x$ head the remaining column and row in that block, respectively.  Equating the $00$ column and row sums in that block gives $a=x$, violating the sudoku condition.
\hspace{\fill} $\Box$
\bigskip

\begin{lemma}
The hypothesis on $\mathrm{exp}(G)$ cannot be omitted from Theorem \ref{ag}.  Furthermore, hypotheses \ref{tppa} and \ref{mspa} cannot both be omitted.
\end{lemma}

\emph{Proof:}  Consider the group $G = \mathbb{Z}_9$ with $k=3$ and $N = \left< 3 \right>$.  It is easy to show that all the hypotheses of Theorem \ref{ag} are satisfied except that $\textrm{exp}(G) = 9$ does not divide $3$.  For example, hypothesis \ref{mspa} is (because the operation is addition mod $9$) $(0 + N) + (1 + N) + (2 +N) = 3 + N = N$.  Therefore, by Lemma \ref{nomcs}, the hypothesis on $\textrm{exp(G)}$ is necessary.

Turn now to the group $G = \mathbb{Z}_2 \times \mathbb{Z}_2$ with $N = \mathbb{Z}_2 \times \{0\}$ (or any subgroup of order 2, for that matter).  Neither of hypotheses \ref{mspa} nor \ref{tppa} are satisfied in this case, but all of the others hold.  By Lemma \ref{nomcs}, those hypotheses cannot both be omitted.
\hspace{\fill} $\Box$  

\section{Other Applications}

In this section, we give two more (of several possible) applications of Theorem \ref{mcs}.  First, we prove every finite group is isomorphic to a subgroup of a finite group having a pandiagonal magic Cayley-sudoku table. Next, we give a class of non-abelian groups of prime power order (i.e. p-groups) having pandiagonal magic Cayley-sudoku tables.

\begin{appl}[Arbitrary Groups] Any finite group is isomorphic to a subgroup of a finite group having a pandiagonal magic Cayley-sudoku table.
\end{appl}

\emph{Proof.}  Let $H$ be any finite group.  Define $G := H \times H \times \mathbb{Z}_{|H|} \times \mathbb{Z}_{|H|}$ and $N := 1 \times 1 \times \mathbb{Z}_{|H|} \times \mathbb{Z}_{|H|}$.  We claim Theorem \ref{mcs} with $k = |H|^2$ applies. Clearly, $N \leq Z(G)$ and $\textrm{exp}(G) = |H|$ divides $k$.  Hypothesis \ref{tpp} holds by Lemma \ref{s2}. Take any ordering $H = [h_1, h_2, \ldots , h_n]$ of the elements of $H$ and define $T$ to be the elements of $H \times H \times 1 \times 1$ ordered lexicographically. That is, $(h_{i_1}, h_{j_1},1,1)$ precedes $(h_{i_2}, h_{j_2},1,1)$ if and only if ${i_1} < {i_2}$ or else ${i_1} = {i_2}$ and ${j_1} < {j_2}$.  In other words, $T := [(h_1, h_1, 1, 1), (h_1 , h_2, 1, 1), \ldots, (h_1, h_n, 1, 1), (h_2, h_1, 1, 1), (h_2, h_2, 1, 1) \ldots, \\ (h_2, h_n, 1, 1) , \ldots, (h_n, h_1, 1, 1), (h_n h_2, 1, 1) \ldots , (h_n, h_n, 1, 1)]$.  It is easy to check $T$ is a complete set of left coset representatives of $N$ in $G$ since $H \times H \times 1 \times 1$ is a complement of $N$ in $G$. To check hypothesis \ref{msp}, fix $1 \leq i, j \leq n$ and note the desired product is $\prod_{\ell =1}^n \left( \prod_{m = 1}^n (h_\ell, h_m, 1, 1)(h_i, h_j, 1, 1) \right) = \prod_{\ell =1}^n \left( (h_\ell h_i)^n, \prod_{m = 1}^n (h_m h_j), 1, 1 \right)$ $ = \prod_{\ell =1}^n \left( 1, \prod_{m = 1}^n (h_m h_j), 1, 1 \right) = \\ \left( 1, (\prod_{m = 1}^n (h_m h_j))^n, 1, 1 \right) = (1,1,1,1)$, where we used the fact that since $n = |H|$, $h^n = 1$ for any $h \in H$.  Thus, $G$ has a pandiagonal magic Cayley-sudoku table and it visibly has a subgroup isomorphic to $H$.

\bigskip
In the following application, the reader unfamiliar with extra-special p-groups will find all the necessary information listed in the first step of the proof.

\begin{appl}[Certain p-Groups]  Suppose $E$ is an extra-special p-group of order $p^3$, where $p$ is an odd prime, and $G = E \times \mathbb{Z}_p$.  Then there is a complete set of left coset representatives, $T$, of $N := Z(G)$ in $G$ such that the table of $G$ based on $N$ and $T$ is a pandiagonal magic Cayley-sudoku table.
\end{appl}

\emph{Proof.} Assume $p$ is an odd prime. In order to construct $T$, some technical results about $E$ are needed, which we number for convenience.

\begin{enumerate}
\item $E$ is either $\left<a, b : a^{p^2} = b^p =1, ab = ba^{1+p}  \right>$, called the exponent $p^2$ case, or $\left< a, b, c : a^p = b^p = c^p = 1, ac = ca, bc = cb, ab = bac \right>$, the exponent $p$ case.  Furthermore, $Z(E) = <a^p>$ or $<c>$, in the exponent $p^2$ and exponent $p$ cases, respectively.

This follows from the classification of extra-special p-groups. See [\ref{gor}] Theorem 5.5.1, for example.

Finally, $E/Z(E)$ is abelian since its order is $p^2$.

\item \label{z} $ab = baz$ for some $z \in Z(E)$ with $|z|=p$.

In the exponent $p^2$ case, $ab = ba^{1+p} = baa^p$ and $a^p \in Z(E)$.  In the exponent $p$ case, $ab = bac$ and $c \in Z(E)$.  In each case, $|z|=p$.
\bigskip

\item \label{lcr}  The set $\{b^j a^i : 0 \leq i, j \leq p-1 \}$ is a complete set of left coset representatives of $Z(E)$ in $E$.

First, suppose $b^j a^i \in Z(E)$ for some $-(p-1) \leq i, j \leq p-1$.  We aim to show $i=j=0$. Using step \ref{z}, we have $b^j a^i = (b^j a^i)^{b^{-1}} = b^j (a^{b^{-1}})^i = b^j (az)^i = b^j a^i z^i$.  It follows that $1 = z^i$, which implies $p = |z|$ divides $i$ and so $i=0$.  Therefore, $b^j \in Z(E) \cap \left< b \right> = 1$, so $p = |b|$ divides $j$.  Thus, $j=0$.

Suppose $a^ib^j Z(E) = a^\ell b^m Z(E)$  where $0 \leq i, j, \ell, m \leq p-1$, then $Z(E) = (a^ib^j)^{-1}(a^\ell b^m) Z(E) = a^{\ell-i} b^{m-j} Z(E)$, the last equality a consequence of $E/Z(E)$ being abelian. In other words, $a^{\ell-i} b^{m-j} \in Z(E)$. From the preceding paragraph, $i = \ell$ and $j=m$.  Consequently, the $p^2$ elements of $\{b^j a^i : 0 \leq i, j \leq p-1 \}$ each represent a different left coset of $Z(E)$ in $E$.  Since $[E:Z(E)]=p^2$, we have a complete set of left coset representatives as was to be shown.

The next three steps are easy to prove using induction.  Each one uses its predecessor.

\item For every integer $i \geq 0$, $a^ib = ba^iz^i$.

\item \label{abz} For all integers $i, j \geq 0$, $a^ib^j = b^ja^i(z^i)^j$.

Fix $i$ and induct on $j$.

\item \label{prod} For all integers $n, i, m \geq 0$, 
\\ $\prod_{\ell = 0}^n (b^{\ell+i} a^m) = \left( \prod_{\ell = 0}^n b^{\ell+i} \right) (a^m)^{n+1} \prod_{\ell = 0}^n (z^m)^{\ell(\ell+i)}$.

Induct on $n$.
\bigskip

The next step is for eventual use in proving hypothesis \ref{msp} of Theorem \ref{mcs}.

\item \label{pmsp} For all integers $i, j \geq 0$, $\prod_{k = 0}^{p-1} \left( \prod_{\ell = 0}^{p-1} (b^\ell a^k)(b^i a^j)  \right) = 1$.

To establish this, we begin by calculating the inner product, first utilizing step \ref{abz} and, lastly, the fact that $z \in Z(E)$ has odd order $p$.
\begin{eqnarray*}
\prod_{\ell = 0}^{p-1} (b^\ell a^k)(b^i a^j) & = & \prod_{\ell = 0}^{p-1} b^\ell (a^k b^i) a^j \\
\null & = & \prod_{\ell = 0}^{p-1} b^{\ell + i} a^k (z^k)^i a^j \\
\null & = & \left( \prod_{\ell = 0}^{p-1} b^{\ell + i} a^{k+j} \right) \prod_{\ell = 0}^{p-1} z^{ki} \\
\null & = & \left( \prod_{\ell = 0}^{p-1} b^{\ell + i} a^{k+j} \right) (z^{ki})^{\frac{(p-1)p}{2}} \\
\null & = & \prod_{\ell = 0}^{p-1} b^{\ell + i} a^{k+j}
\end{eqnarray*}

Further simplification follows from step \ref{prod}, temporarily setting $m:=k+j$ for notational convenience.  Recall that $b$ and $z$ each have odd order $p$.
\begin{eqnarray*}
\prod_{\ell = 0}^{p-1} b^{\ell + i} a^{k+j} & = & \left( \prod_{\ell = 0}^{p-1} b^{\ell+i} \right) (a^m)^p \prod_{\ell = 0}^{p-1} (z^m)^{\ell(\ell+i)} \\
\null & = & b^{\frac{(p-1)p}{2} + p i } a^{pm} (z^m)^{\frac{(p-1)p(2p-1)}{6} + \frac{(p-1)p}{2}i } \\
\null & = & a^{pm} z^{\frac{(p-1)p(2p-1)}{6}} \\
\null & = & \left\{  \begin{array}{ll}
                    a^{pm} (z^m)^5 = a^{pm} (z^m)^2 & p = 3 \\
                    a^{pm} & p \geq 5
                    \end{array}   \right.
\end{eqnarray*}

Thus, the desired double product becomes
\begin{eqnarray*}
\prod_{k = 0}^{p-1} \left( \prod_{\ell = 0}^{p-1} (b^\ell a^k)(b^i a^j) \right) & = & \left\{  \begin{array}{ll}
                    \prod_{k = 0}^{p-1} a^{pm} (z^m)^2 & p = 3 \\
                    \prod_{k = 0}^{p-1} a^{pm} & p \geq 5
                    \end{array}   \right.
\end{eqnarray*}

We split into cases.  First consider the exponent $p^2$ case where $z = a^p$ (see step \ref{z}).

When $p=3$, $a^{pm} (z^m)^2 = (z^m)^3 = 1$ and we are finished.  When $p \geq 5$, the double product is $\prod_{k = 0}^{p-1} z^m = \prod_{k = 0}^{p-1} z^{k+j}  = z^{\frac{(p-1)p}{2} + pj} = 1$.

Turn next to the exponent $p$ case, where $a^p = 1$.  When $p \geq 5$ the double product is clearly $1$.  When $p = 3$, we have $\prod_{k = 0}^{p-1} a^{pm} (z^m)^2 = \prod_{k = 0}^{p-1} (z^{k+j})^2 = z^{(p-1)p + 2pj} = 1$.
\end{enumerate}

We will now show Theorem \ref{mcs} applies with $N =  Z(G) = Z(E) \times \mathbb{Z}_p$ and $k = p^2 = |N|$. First note that $\textrm{exp}(G) = \textrm{exp}(E)$ is either $p$ or $p^2$. The first two hypotheses of the theorem are now clear. Hypothesis \ref{tpp} holds by Lemma \ref{s2}.  Let $T$ be the list $[(b^\ell a^k,1)$ : $0 \leq \ell, k \leq p-1]$ ordered lexicographically on $\ell$ and $k$, which is to say $(b^{\ell_1} a^{k_1},1)$ precedes $(b^{\ell_2} a^{k_1},2)$ if and only if $\ell_1 < \ell_2$ or else $\ell_1 = \ell_2$ and $k_1 < k_2$.  That $T$ is a complete set of left coset representatives of $N$ in $G$ follows easily from step \ref{lcr} and the fact that $N = Z(E) \times \mathbb{Z}_p$.  Furthermore, $T$ satisfies Hypothesis \ref{msp} by step \ref{pmsp}, which completes the proof.
\hspace{\fill} $\Box$

\section*{References}

\begin{enumerate}
\item \label{csw}  J. Carmichael, K. Schloeman, and M. B. Ward, Cosets and Cayley-Sudoku Tables, Math. Mag., \textbf{83} (2010) 130-139.  \\ DOI: 10.4169/002557010X482899.
\item \label{gor} D. Gorenstein, Finite Groups, Harper \& Row, 1968. 
\item \label{LW} J. Lorch and E. Weld, Modular Magic Sudoku, Involve 5 (2012) 173-186.  DOI: 10.2140/involve.2012.5.173.
\item \label{P} L. J. Paige, A Note on Finite Abelian Groups, Bull. Amer. Math. Soc. 53 (1947) 590-593. \\
 \verb+https://projecteuclid.org/euclid.bams/1183510800+

\end{enumerate}

\end{document}